# QUASIPOTENTIALS WITH MORE THAN TWO VARIABLES: NEW EVALUATION AT EQUILIBRIUM POINTS OF THE DRIFT

Dietrich Ryter

RyterDM@gawnet.ch

Midartweg 3   CH-4500 Solothurn  Switzerland

Phone   +4132 621 13 07

The relevant quasipotential near an equilibrium point is determined by a new linear matrix equation, with less unknowns than an existing nonlinear one. This also assures the asymptotic fulfillment of the Fokker-Planck equation, even globally due to the second term in the noise strength.
An auxiliary result for the exit problem is derived as well.

*Key words:*  Fokker-Planck equation; weak noise; quasipotential; exit problem

# I. Introduction

The recent paper [1] presented a new evaluation method for the quasipotential (QP) with two variables. The full extension to more variables seems unfeasible, but some important findings can be generalized; this mainly concerns the local approximation at equilibrium points (EP) of the drift (which can be continued to the remaining space by the usual Hamiltonian system). The new approach yields a linear matrix equation which is superior to the well-known algebraic Riccati equation, because

(1) it is always linear

(2) it is more general than the *linear* Ricatti equation

(3) it involves less parameters, since the unknown matrix is antisymmetric rather than symmetric.

This approach solves the Fokker-Planck equation (FPE) at the EP, and thereby also guarantees the local regularity of the function which extends the QP to the next order in the noise strength and entails the global asymptotic fulfillment of the FPE. An example shows that extra solutions of the nonlinear Riccati equation are irrelevant.

These findings are of immediate use for constructing the Gaussian pdf near a point attractor, but also important for the exit problem from a basin of attraction, when the separatrix has an attractive point : by [1] the "associated drift" is crucial here, and its local approximation is obtained as well, including its unstable direction, on which numerical integrations are conveniently started.

# II. Background

A stationary FPE [2,3] with $n$ variables $x^i$, $i = 1,...,n$, drift components $a^i(\vec{x})$, and with the diffusion matrix $2\varepsilon \underline{D}(\vec{x})$ (symmetric and nonnegative) can be written as

$$\nabla \cdot (-\vec{a}w + \varepsilon \underline{D}\nabla w) = \rho w - \vec{a}\cdot\nabla w + \varepsilon \nabla\cdot(\underline{D}\nabla w) = 0 , \qquad (2.1)$$

where $\rho(\vec{x}) := -\nabla\cdot\vec{a}$ means the contraction of $\vec{a}$, and $\varepsilon$ exhibits the noise strength.

The noise-induced drift $\varepsilon\,\partial D^{ij}/\partial x^j$ (with summation over $j$) has been included in (2.1), see e.g. [2,4].

The drift $\vec{a}$ is further assumed to have an EP with a linear decay, first in the sense that the matrix $\underline{M}$, consisting of the local row vectors $\nabla a^i$, is not singular.

When the solution $w$ is written as

$$w(\vec{x}) = N\exp\{-\varepsilon^{-1}[\phi(\vec{x}) + \varepsilon\varphi(\vec{x})]\} = N\exp[-\varphi(\vec{x})]\exp[-\phi(\vec{x})/\varepsilon] \qquad (2.2)$$

[5,6], the FPE (2.1) becomes

$$\varepsilon^{-1}(\vec{a} + \underline{D}\nabla\phi)\cdot\nabla\phi + r - \tilde{\vec{a}}\cdot\nabla\varphi + \varepsilon[\nabla\varphi\cdot(\underline{D}\nabla\varphi) - \nabla\cdot(\underline{D}\nabla\varphi)] = 0 \qquad (2.3)$$

where

$$r := \rho - \nabla\cdot(\underline{D}\nabla\phi) \qquad (2.4)$$

$$\tilde{\vec{a}} := -(\vec{a} + 2\underline{D}\nabla\phi) \qquad (2.5)$$

($\tilde{\vec{a}}$ was called "associated drift" in [1]).

For small $\varepsilon$ it is natural to determine the QP or "eikonal" $\phi$ by the Freidlin equation

$$(\vec{a} + \underline{D}\nabla\phi)\cdot\nabla\phi = 0 , \qquad (2.6)$$

which is of the first order, but quadratic in the derivatives of $\phi$.

Except for "exact" solutions with $r \equiv 0$ (which include the cases with detailed balance [3]) $\varphi$ is not a constant, but when it satisfies

$$\tilde{\vec{a}}\cdot\nabla\varphi = r \qquad (2.7)$$

[6], the FPE is fulfilled to $O(\varepsilon)$. At an EP $\tilde{\vec{a}} = \vec{0}$. This means that $\varphi$ is only regular, when $r$ vanishes as well. It will be shown that this holds indeed for the relevant QP, but not in general with further solutions of (2.6).

The usual way of solving (2.6) is to consider the Hamiltonian $H = p_i(a^i + D^{ij}p_j)$

with the momenta $p_i := \partial \phi / \partial x^i$, and to integrate

$$\dot{x}^i = \partial H / \partial p_i = a^i + 2D^{ij} p_j \qquad (= -\tilde{\vec{a}}) \tag{2.8}$$

$$\dot{p}_i = -\partial H / \partial x^i = -p_k (\partial a^k / \partial x^i + p_j \partial D^{jk} / \partial x^i) \tag{2.9}$$

as well as $\dot{\phi} = p_i \dot{x}^i$,

Starting conditions near an EP are provided by local analytical solutions. We mention that the second derivatives of $\phi$ (required in (2.7) via (2.4)) can be computed together with the integration of (2.8), (2.9), by the method of [7].

With two variables ($n = 2$) the Hamilton method can be substituted by the approach of [1], which avoids the doubling of the variables.

### III. An Alternative Evaluation of the QP

Actually (2.6) states that the "conservative drift" $(\vec{a} + \underline{D}\nabla\phi) := \vec{a}_c$ is orthogonal to $\nabla\phi$, so that $\vec{a}_c = \underline{A}\nabla\phi$, with an antisymmetric $\underline{A}(\vec{x})$. This entails

$$\vec{a} = (-\underline{D} + \underline{A})\nabla\phi . \tag{3.1}$$

As a consequence

$$\nabla\phi = (-\underline{D} + \underline{A})^{-1}\vec{a} , \tag{3.2}$$

and the remaining problem is to determine $\underline{A}(\vec{x})$.

While $\vec{a}$ of (3.1) satisfies (2.6) for any $\underline{A}(\vec{x})$, it is obvious that not every solution of (2.6) agrees with (3.1): a counterexample is $\nabla\phi \equiv \underline{0}$. A necessary condition for $\phi$ at an EP can be obtained by taking all first derivatives of (3.2). With the second derivatives of $\phi$ arranged as a matrix $\underline{S}$, it follows that

$$\underline{S} = (-\underline{D} + \underline{A})^{-1}\underline{M} , \tag{3.3}$$

recall that the $i$-th row of $\underline{M}$ is $\nabla a^i$. The rank of $\underline{S}$ is thus that of $\underline{M}$.

The complete $\underline{A}(\vec{x})$ was evaluated in [1] for the case of two variables ($x^1 := x$, $x^2 := y$), where $\underline{A}$ is given by a single function $\chi(x,y)$

$$\underline{A} = \chi \begin{pmatrix} 0 & 1 \\ -1 & 0 \end{pmatrix} . \tag{3.4}$$

The request that the $\phi_x, \phi_y$ (as given by (3.2)) must satisfy the gradient condition $(\phi_x)_y = (\phi_y)_x$ resulted in a *quasilinear* pde of the first order for $\chi$, with characteristics in the plane, given by the associated drift $\tilde{\vec{a}}$.

A similar procedure for $n > 2$ would result in a system of coupled partial pde's, which does not look promising. We rather focus on $\underline{A}$ at EP's (in the next Chapter), and merely observe the following:

1) Taking the divergence of (3.1) shows that $\rho + \nabla \cdot (\underline{D} \nabla \phi) = r$ equals a term which vanishes when $\nabla \phi = \vec{0}$ or when $\underline{A}$ is constant (mind that the trace of a product of $\underline{A}$ with a symmetric matrix is zero). At an EP ($\nabla \phi = \vec{0}$) $r$ vanishes thus indeed, so that $\varphi$ is regular there (apparently a new finding for $n > 2$, but only valid with (3.3)).

2) $\underline{A}$ cannot be finite on a possible *limit cycle* of $\vec{a}$, since $\nabla \phi$ vanishes there (by the Lyapunov property of $\phi$), but not $\vec{a}_c$, i.e. the drift on the cycle. The above argument does thus not apply and $r \neq 0$ there. Since also $\tilde{\vec{a}} = -\vec{a}_c \neq \vec{0}$, $\varphi$ is well defined, and $\exp(-\varphi)$ is essentially given by $|\vec{a}_c|^{-1}$. We mention without a proof that for $n = 2$ there is a simple relation between the gradient of $\chi^{-1}$ and the second derivative of $\phi$ normal to the cycle.

### IV. $\underline{A}$ at an Equilibrium Point

4.1 *The equation for $\underline{A}$*

The matrix $\underline{S}$ is symmetric. It is remarkable that this already determines $\underline{A}$ in (3.3), when both $\underline{S}$ and $\underline{M}$ are non-singular (as supposed for the present). To see this, consider the inverse of (3.3) $\underline{S}^{-1} = \underline{M}^{-1}(-\underline{D} + \underline{A})$, which must equal its transpose $(-\underline{D} + \underline{A})^T (\underline{M}^{-1})^T = -(\underline{D} + \underline{A})(\underline{M}^{-1})^T$, so that

$$\underline{M}^{-1}\underline{A} + \underline{A}\,(\underline{M}^{-1})^T = \underline{M}^{-1}\underline{D} - \underline{D}(\underline{M}^{-1})^T \ .$$

Multiplication from the left by $\underline{M}$ and from the right by $\underline{M}^T$ finally yields

$$\underline{A}\,\underline{M}^T + \underline{M}\,\underline{A} = \underline{D}\,\underline{M}^T - \underline{M}\,\underline{D} \ . \tag{4.1}$$

It is easily seen that both sides are antisymmetric. This relation is in fact a *linear* system of $n(n-1)/2$ equations for the $n(n-1)/2$ independent elements of $\underline{A}$. The fact that $\underline{M}^{-1}$ does not occur in (4.1) may even admit singular $\underline{M}$, when the solution of (4.1) is unique and when the resulting $(-\underline{D} + \underline{A})^{-1}$ exists (which is not evident when $\underline{D}$ is singular). An example for this is the Kramers model with a flat bottom or barrier: both $\underline{M}$ and $\underline{S}$ are then singular, but $\underline{A}$, as well as $(-\underline{D} + \underline{A})^{-1}$, are well defined and correct, see below.

Note that $\underline{A}$ vanishes when $\underline{M}$ is symmetric and commutes with $\underline{D}$.

One can also obtain an equation for $\underline{S}^{-1}$ which avoids $\underline{A}$. To this end consider again the inverse of (3.3) $\underline{S}^{-1} = \underline{M}^{-1}(-\underline{D} + \underline{A})$ and multiply this from the left by $\underline{M}$, which yields $\underline{M}\,\underline{S}^{-1} = -\underline{D} + \underline{A}$. Addition of the transpose $\underline{S}^{-1}\underline{M}^T = -\underline{D} - \underline{A}$ results in

$$\underline{M}\,\underline{S}^{-1} + \underline{S}^{-1}\underline{M}^T = -2\underline{D} \ . \tag{4.2}$$

This is again a linear matrix equation, now for the symmetric $\underline{S}^{-1}$. It involves $n$ more independent elements than (4.1), and explicitly requires the existence of $\underline{S}^{-1}$.

In [7] it was shown how the second derivatives of $\phi$ develop along the bicharacteristics of the Hamiltonian system (the definition of $\underline{D}$ differs by a factor 2).

At an EP this resulted in the nonlinear equation

$$\underline{S}\,\underline{M} + \underline{M}^T\underline{S} + 2\underline{S}\,\underline{D}\,\underline{S} = \underline{0} \qquad (4.3)$$

of the algebraic Riccati type. When $\underline{S}^{-1}$ exists, this reduces to (4.2). Further solutions of (4.3) have a lower rank, as for example the trivial $\underline{S} = \underline{0}$. For the Kramers model examples with rank 1 (see below) have neither a physical nor a probabilistic meaning.

4.2 *Solving (4.1) for $\underline{A}$*

a) $n = 2$

In terms of $x^1 := x$, $x^2 := y$; $a^1 := a$, $a^2 := b$  $\underline{M}$ reads

$$\begin{pmatrix} a_x & a_y \\ b_x & b_y \end{pmatrix}.$$

With (3.4) the lefthand side of (4.1) becomes

$$(a_x + b_y)\,\chi \begin{pmatrix} 0 & 1 \\ -1 & 0 \end{pmatrix},$$

so that $\chi$ is determined whenever $a_x + b_y \neq 0$ (i.e. $\rho \neq 0$; equality is possible at a hyperbolic point). The result (5.8) of [1] is recovered, in particular $\chi = (b_x - a_y)/(a_x + b_y)$ for $\underline{D} = \underline{I}$, and $\chi = 1$ in the Kramers model (where $\rho = \gamma > 0$), irrespective of $U''$, thus also for $U'' = 0$ where both $\underline{M}$ and $\underline{S}$ are singular.

b) $n \geq 3$

A solution in a closed form is presently not known. For the numerical evaluation consider the antisymmetric matrices $\underline{e}_{ik}$ with the elements 1 at $i < k$ and -1 with $i, k$ interchanged, and with zeros elsewhere. Clearly they are a basis in the space of the antisymmetric matrices. The aim is to evaluate the coefficients $\alpha_{ik}$ in

$$\underline{A} = \sum_{i<k} \alpha_{ik}\,\underline{e}_{ik}.$$

Inserting this into (4.1) - and representing the righthand side accordingly - yields the required linear equations by annihilating the resulting prefactors of all $\underline{e}_{ik}$.

The case $\rho (= -tr\underline{M}) = 0$ seems still problematic, since $\varepsilon$ does not occur in the corresponding local FPE $\nabla \cdot (\underline{D} \nabla w) = 0$.

## V. An Auxiliary Result for the Exit Problem

The exit problem, as treated in [1], requires the knowledge of the analogue of $\underline{M}$ for the drift with the reversed conservative part, i.e. of

$$\tilde{\vec{a}} = (-\underline{D} - \underline{A}) \nabla \phi$$

(with $\phi$ and thus $\underline{S}$ unchanged). Clearly by (3.3) $\underline{M} = (-\underline{D} + \underline{A})\underline{S}$ and therefore $\underline{\tilde{M}} = (-\underline{D} - \underline{A})\underline{S}$, so that

$$\underline{\tilde{M}} = \underline{M} - 2\underline{A}\underline{S} \ . \tag{5.1}$$

This is to be used at an EP on the separatrix of $\vec{a}$, which is attractive within the separatrix and repulsive across it; accordingly one eigenvalue of $\underline{M}$ is realvalued and positive, while the real parts of the others are negative. To obtain the eigenvalues of $\underline{\tilde{M}}$ one may solve $\underline{M} = (-\underline{D} + \underline{A})\underline{S}$ for $\underline{A}$ and use the antisymmetry, which yields $-\underline{A} = \underline{D} + \underline{S}^{-1} \underline{M}^T$. Inserted into $\underline{\tilde{M}} = (-\underline{D} - \underline{A})\underline{S}$ this results in

$$\underline{\tilde{M}} = \underline{S}^{-1} \underline{M}^T \underline{S} \ . \tag{5.2}$$

The eigenvalues of $\underline{\tilde{M}}$ are thus those of $\underline{M}$. When $\underline{S}$ is singular, the traces of $\underline{\tilde{M}}$ and $\underline{M}$ are still the same (by $tr(\underline{A}\underline{S}) = 0$), and also the determinants, since $\underline{\tilde{M}} = (-\underline{D} - \underline{A})(-\underline{D} + \underline{A})^{-1} \underline{M}$ and $(\underline{D} - \underline{A}) = (\underline{D} + \underline{A})^T$. The unstable eigendirection of $\underline{\tilde{M}}$ (relevant for starting the numerical evaluation of $\phi$, since (2.8) amounts to $\dot{\vec{x}} = -\tilde{\vec{a}}$) is therefore given by $\underline{S}^{-1} f$, where $f$ is the eigenvector of $\underline{M}^T$ with the positive eigenvalue $\lambda_+$. By $\underline{M}^T f = \lambda_+ f$ and $\underline{S}^{-1}\underline{M}^T = -\underline{D} - \underline{A}$ (see before (4.2))
this has the direction of $-(\underline{D} + \underline{A})f$, which does not involve $\underline{S}^{-1}$.

## VI  The Kramers Example

### 6.1 *The model*

This model [8,9] with $n = 2$ describes a massive particle moving in a potential $U(x)$.
With unit mass and with temperature $\varepsilon$ the equation of motion is

$$\dot{v} = -\gamma v - U'(x) + (2\varepsilon\gamma)^{1/2}\xi \qquad (\xi \text{ being standard white noise}).$$

In the variables $x, v$ the drift is $a^1 = v$, $a^2 = -\gamma v - U'(x)$, while $D^{22} = \gamma$ and the other elements of $\underline{D}$ vanish. The equation (2.6) for $\phi$ reads

$$v\phi_x - (\gamma v + U')\phi_v + \gamma(\phi_v)^2 = 0 \tag{6.1}$$

and is solved by the well-known $\phi_{eq} = U(x) + v^2/2$. The EP's are given by $U' = 0 = v$, with the corresponding

$$\underline{M} = \begin{pmatrix} 0 & 1 \\ -U'' & -\gamma \end{pmatrix},$$

while for $\phi_{eq}$ $\underline{S}$ is diagonal with $S_{11} = U''$ and $S_{22} = 1$.

### 6.2 *A degenerate quasipotential*

In the quadratic approximation at the bottom ($U'' \geq 0$) or at a barrier ($U'' \leq 0$) of $U$ (with $x = v = 0 = \phi$ there) two further solutions of (6.1) are

$$\phi_{\pm} = [U''(1-\beta)x^2 + 2U''\gamma^{-1}xv + \beta v^2]/2 = (2\beta)^{-1}(U''\gamma^{-1}x + \beta v)^2 \tag{6.2}$$

$$\text{where} \quad 2\beta_{\pm} = 1 \pm (1 - 4U''\gamma^{-2})^{1/2}.$$

Both $\phi_{\pm}$ always exist at a barrier, and also at the bottom when the local oscillation is overdamped; in the latter case (6.2) holds globally when $U''(x)$ is constant. The respective $w$ (without $\varphi$) is however not concentrated at the bottom (!) and cannot be normalized. The $\underline{S}$ corresponding to (6.2) has rank 1 and solves (4.3); it yields

$$r = \rho - \nabla \cdot (\underline{D}\nabla\phi) = \gamma(1-\beta_\pm) \neq 0 \quad \text{(unless } U''=0 \text{ where } \beta_\pm = 1\text{)}.$$

As a consequence $\varphi$ is singular at the EP. [The same holds for $\phi_0 \equiv 0$ (rank 0), where $r = \gamma > 0$].

6.3 *Failure of the Minimum Principle*

Since $\phi$ is the action function of a Hamiltonian, it is worthwhile to check whether the minimum principle selects the relevant solution of (2.6) near an EP, namely the one with the maximum rank of $\underline{S}$ ($\phi_{eq}$ in the Kramers case). The answer is clearly negative, since $\phi_{eq} \geq U_{\min}$, and $\phi_0 = const. = U_{\min}$ would be smaller; also $\phi_\pm$ (when they exist) are smaller than $\phi_{eq}$ on one side of the line through $(x_{\min},0)$ with the slope $\beta_\pm\gamma$. Near the barrier a rigorous application of the minimum principle would even entail a complicated "patchwork" solution (continuous, but with a discontinuous gradient at the seems), since the minimum $\phi_{eq}$, $\phi_\pm$, $\phi_0$ is not the same in different directions from the EP. The nondifferentiable QP's obtained by this principle [10-12] might therefore not be relevant for the FPE.